\documentclass[]{amsart}

\usepackage{amssymb}
\usepackage{amsmath}
\usepackage{indentfirst}
\usepackage{amscd}
\usepackage{graphicx}

\newtheorem{thm}{Theorem}[section]
\newtheorem{lem}[thm]{Lemma}
\newtheorem{prop}[thm]{Proposition}

\newtheorem{cor}[thm]{Corollary}

\newtheorem{rem}[thm]{Remark}

\newtheorem{prob}[thm]{Problem}

\newcommand{\R}{\ensuremath{\mathbb{R}}}

\newcommand{\Z}{\ensuremath{\mathbb{Z}}}

\newcommand{\Hom}[2]
{ \ensuremath{#1^{ \mathrm{#2} }}  }
\newcommand{\Coh}[2]
{ \ensuremath{#1_{ \mathrm{#2} }}  }

\newcommand{\bou}{\ensuremath{\partial}}
\newcommand{\cob}{\ensuremath{\delta}}
\begin{document}

%%%%%%%%%%%%%%%%%%%%%%%%%
% Subject classification 
%%%%%%%%%%%%%%%%%%%%%%%%%

\date{December 5, 2005}

%%%%%%%%%
% Title
%%%%%%%%%

\title[]
{Inequivalent surface-knots with 
the same knot quandle}

%%%%%%%%%%%%%%%%%%%%%%%%%%%%%%
% Author names and addresses 
%%%%%%%%%%%%%%%%%%%%%%%%%%%%%%

\author{Kokoro Tanaka}
\address{Graduate School of Mathematical Sciences, 
University of Tokyo, 3-8-1 Komaba Meguro, 
Tokyo 153-8914, Japan}
\email{k-tanaka@ms.u-tokyo.ac.jp}

%%%%%%%%%%%%%
% Dedication
%%%%%%%%%%%%%

%\dedicatory{}

%%%%%%%%%%%%%
% Abstract 
%%%%%%%%%%%%%
\begin{abstract}
We have a knot quandle and a fundamental class
as invariants for a surface-knot.
These invariants can be defined 
for a classical knot in a similar way,
and it is known that the pair of them is
a complete invariant for classical knots.
In this paper, we compare a situation in surface-knot theory 
with that in classical knot theory, and prove the following:
There exist arbitrarily many inequivalent 
surface-knots of genus $g$ with the same knot quandle,
and there exist two inequivalent surface-knots of genus $g$
with the same knot quandle and with the same fundamental class.
\end{abstract}

\maketitle

%%%%%%%%%%%%%%%%%%%%%%%%%%%%%%%%%%%%%%%%%%%%%%%%%%%%%%%%%%%%%%%%%%%%%%%%%
% end Topmatter
%%%%%%%%%%%%%%%%%%%%%%%%%%%%%%%%%%%%%%%%%%%%%%%%%%%%%%%%%%%%%%%%%%%%%%%%%

%%%%%%%%%%%%%%%%%%%%%%%%%%%%%%%%%%%%%%%%%%%%%%%%%%%%%%%%%%%%%%%%%%%%%%%%%
% body of paper
%%%%%%%%%%%%%%%%%%%%%%%%%%%%%%%%%%%%%%%%%%%%%%%%%%%%%%%%%%%%%%%%%%%%%%%%%
\section{Introduction}
We consider a {\it knot quandle} \cite{Joy, Mat}, $Q(F)$, 
and a {\it fundamental class} \cite{CKS} (cf. \cite{Tana2}), 
$[F] \in \Hom{H}{Q}_3(Q(F))$,  
as invariants of a surface-knot $F$,
where a surface-knot means an oriented closed connected surface
embedded in $\R^4$.
The fundamental class can be considered as a universal object concerning to
a {\it quandle cocycle invariant} (See Section~\ref{subsec-qc_inv}).
When the invariants are given, what we want to know might be 
the following:
\begin{itemize}
\item
What kind of information can be extracted from them?
\item
How powerful are they?
\end{itemize}

For the first question, 
it is known in \cite{Joy, Mat}
that the knot quandle of a surface-knot $F$
can recover information of the knot group $\pi_1(\R^4 \setminus F)$,
for example.
There are some relation of the knot quandle
to the braid index \cite{Tana},
to the unknotting number \cite{Iwa2} and 
to the sheet number \cite{Sai-Sat}.
There are also some relation of the fundamental class 
to the non-invertibility \cite{CJKLS, AS, Iwa},
to the triple point number \cite{Sa-Shi, Sa-Shi2, Hat, Tana2},
to the triple point cancelling number \cite{Iwa2}, and
to the ribbon concordance \cite{CSS}.

For the second question, it is known in \cite{BL} that
the knot quandle can distinguish all elements of a class
of twist-spun $S^2$-knots obtained from torus knots,
for example.
In this paper, we focus on the second question
and compare a situation in surface-knot theory with 
that in classical knot theory.

%%%%%%%%%%%%%%%%%%%
\subsection{The case of classical knots}\label{subsec-classical}
Similarly, we have a knot quandle $Q(k)$ and 
a fundamental class $[k] \in \Hom{H}{Q}_2(Q(k))$ as invariants
of a classical knot $k$ (cf. \cite{Eis}). 
For a classical knot $k$,
let $-k$ denote the classical knot obtained from $k$ by 
reversing the orientation,
and $k^{\ast}$ denote the mirror image of $k$.
Then the following three facts are known.

\begin{itemize}
\item
\underline{Fact (cf. \cite[Proof of Theorem $9.1$]{CJKS})}:
For a classical knot $k$,
there exists a canonical quandle isomorphism
$\phi : Q(k) \rightarrow Q(-k^{\ast})$ such that
the induced homomorphism
$\phi_{\ast} : \Hom{H}{Q}_2(Q(k)) \rightarrow \Hom{H}{Q}_2(Q(-k^{\ast}))$
satisfies the condition $\phi_\ast [k]=-[-k^{\ast}]$.
\item
\underline{Theorem due to Joyce \cite{Joy} and Matveev \cite{Mat}}:
For classical knots $k$ and $k'$,
if there exists a quandle isomorphism
$\phi : Q(k) \rightarrow Q(k')$,
then $k$ is equivalent to $k'$ or $-(k')^{\ast}$.
\item
\underline{Theorem due to Eisermann \cite{Eis}}:
For classical knots $k$ and $k'$,
if there exists a quandle isomorphism
$\phi : Q(k) \rightarrow Q(k')$
such that the induced homomorphism $\phi_{\ast}$
satisfies the condition $\phi_\ast [k]=[k']$,
then $k$ is equivalent to $k'$.
\end{itemize}

Roughly speaking, Joyce--Matveev's theorem says that
the knot quandle is an almost complete invariant for
classical knots,
and Eisermann's theorem says that the pair of the knot
quandle and the fundamental class is a complete
invariant for them.

\begin{rem}\rm
Eisermann \cite{Eis} also proved:
\begin{itemize}
\item
For a trivial classical knot $k$, we have $\Hom{H}{Q}_2(Q(k)) \cong 0$.
\item 
For a non-trivial classical knot $k$, we have $\Hom{H}{Q}_2(Q(k)) \cong \Z$
and the fundamental class $[k]$ is a generator.
\end{itemize}
On the other hand,
as far as the author knows, there is not so much result about
the structure of $\Hom{H}{Q}_3(Q(F))$ for a surface-knot $F$.
\end{rem}

%%%%%%%%%%%%%%%%
\subsection{Problem setting}
For a surface-knot $F$,
let $-F$ denote the surface-knot obtained from $F$ by 
reversing the orientation,
and $F^{\ast}$ denote the mirror image of $F$.
It is known that the assertion corresponding to the first fact 
in Section~\ref{subsec-classical} 
also holds for a surface-knot $F$, that is,
there exists a canonical quandle isomorphism
$\phi : Q(F) \rightarrow Q(-F^{\ast})$ such that
the induced homomorphism
$\phi_{\ast} : \Hom{H}{Q}_3(Q(F)) \rightarrow \Hom{H}{Q}_3(Q(-F^{\ast}))$
satisfies the condition $\phi_\ast [F]=-[-F^{\ast}]$ 
(cf. \cite[Proof of Theorem $9.2$]{CJKS}).
Then we consider the following problem.

\begin{prob}\label{prob-pre}\rm \hspace{0pt}
\begin{enumerate}
\item[(I)]
Does the assertion corresponding to Joyce--Matveev's theorem 
hold for surface-knots?
\item[(II)] 
Does the assertion corresponding to Eisermann's theorem hold for surface-knots?
\end{enumerate}
\end{prob}

Since the knot quandle does not have information of 
the genus of a surface-knot,
we fix a non-negative integer $g$
and consider the above problem for surface-knots of genus $g$.
To make the problem concrete, 
we consider the following five conditions
for two surface-knots, $F$ and $F'$, of genus $g$:
\begin{itemize}
\item[(i)]
There exists a quandle isomorphism $\phi :Q(F) \rightarrow Q(F')$.
\item[(ii)]
There exists a quandle isomorphism $\phi :Q(F) \rightarrow Q(F')$
such that $$\phi_{\ast}[F]=[F'] \in \Hom{H}{Q}_3(Q(F')).$$
\item[(ii')]
There exists a quandle isomorphism $\phi :Q(F) \rightarrow Q(F')$
such that $$\phi_{\ast}[F]=\pm [F'] \in \Hom{H}{Q}_3(Q(F')).$$
\item[(iii)]
The surface-knot $F$ is equivalent to $F'$.
\item[(iii')]
The surface-knot $F$ is equivalent to $F'$ or $-(F')^{\ast}$.
\end{itemize} 
By definition,
we have (iii) $\Rightarrow$ (ii) $\Rightarrow$ (i),
(ii) $\Rightarrow$ (ii'), and (iii) $\Rightarrow$ (iii').
As mentioned above, we also have (iii') $\Rightarrow$ 
(ii') $\Rightarrow$ (i).
Then we can reformurate Problem~\ref{prob-pre} as follows:

\begin{prob}\label{prob-main}\rm
(Reformultation of Problem~\ref{prob-pre})
\begin{itemize}
\item[(I)]
Does the condition (i) imply the condition (iii')? 
\item[(II)]
Does the condition (ii) imply the condition (iii)?
\end{itemize}
Moreover, by the fact that (iii') $\Rightarrow$ (ii') $\Rightarrow$ (i),
we can divide (I) into two parts.
\begin{itemize}
\item[$(\mathrm{I}_1)$]
Does the condition (i) imply the condition (ii')?
\item[$(\mathrm{I}_2)$]
Does the condition (ii') imply the condition (iii')?
\end{itemize}

\end{prob}

The main result of this paper is to give negative answers 
to Problem~\ref{prob-main}.

\begin{thm}\label{th-quandle}
For a non-negative integer $g$,
there exist arbitrarily many surface-knots of genus $g$ 
such that any two of them satisfy the condition {\rm (i)}
but do not satisfy the condition {\rm (ii')}.
\end{thm}

\begin{thm}\label{th-fund}
For a non-negative integer $g$,
there exist two surface-knots of genus $g$
such that they satisfy the condition {\rm (ii)}
but do not satisfy the condition {\rm (iii')}. 
Moreover, infinitely many such pairs exist.
\end{thm}

Theorem~\ref{th-quandle} gives a negative answer to 
Problem~\ref{prob-main} $(\mathrm{I}_1)$,
and Theorem~\ref{th-fund} gives a negative answer to 
Problem~\ref{prob-main} $(\mathrm{I}_2)$ and (II).

\begin{rem}\rm
It follows from Theorem~\ref{th-quandle} that
there exist arbitrarily many inequivalent 
surface-knots of genus $g$ with the same knot group.
We note that the more stronger assertion is known for 
surface-knots of genus zero:   
There exist infinitely many $S^2$-knot
with the same knot group \cite{Suc}.
\end{rem}

The rest of this paper is organized as follows. 
In Section~\ref{sec-lemma}, 
we review the basic definitions
including knot quandles and fundamental classes of surface-knots,
and give Lemma~\ref{lem-quandle} and Corollary~\ref{cor-fund}
which are keys to proving theorems.
Section~\ref{sec-quandle} and Section~\ref{sec-fund} 
are devoted to proving Theorem~\ref{th-quandle}
and Theorem~\ref{th-fund} respectively.

%%%%%%%%%%%%%%%%%%%%%%%%%%%%%%%%%%%%%
\section{Definitions and Lemmas}\label{sec-lemma}

\subsection{Surface-knots and diagrams}
A {\it surface-knot} is a closed connected oriented 
surface embedded locally flatly in $\R^4$
(or in the $4$-sphere $S^4$).
Two surface-knots are said to be {\it equivalent}
if they are related by an ambient isotopy of $\R^4$.
For a fixed projection $\pi : \R^4 \rightarrow \R^3$,  
by perturbing a surface-link $F$ if necessary,
we may assume that the projection $\pi|_F$ is {\it generic}, that is,
$\pi|_F$ has double points, isolated triple points
and isolated branch points in the image as its singularities. 
A {\it diagram} of a surface-knot is a generic projection image 
equipped with height information,
where one of two sheets along each double point curves
is broken depending on the relative height.
A diagram consists of a collection of sheets, 
and is regarded as a compact oriented surface in $\R^3$.
We refer to \cite{CS-book} for more details.

\subsection{Quandles and knot quandles}
A {\it quandle} \cite{Joy, Mat}, $X$, is a non-empty set 
with a binary operation $(a,b) \rightarrow a*b$ satisfying 
the following axioms.
\begin{enumerate} 
\item[(Q1)]
For any $a \in X$, $a*a=a$.
\item[(Q2)]
For any $a,b \in X$, there is a unique $c \in X$
such that $c*b=a$.
\item[(Q3)]
For any $a,b,c \in X$, we have $(a*b)*c=(a*c)*(b*c)$.
\end{enumerate}
A function $f : X \rightarrow Y$ between quandles is a {\it homomorphism}
if $f(a*b)=f(a)*f(b)$ for any $a,b \in X$.

Let $D$ be a diagram of a surface-link $F$, and 
let $E=\{s_1,\dots,s_m\}$ be the set of all sheets of $D$.
Using the orientation of $F$ and that of $\R^3$,
we give a normal vector to each sheet.
The {\it knot quandle} \cite{Joy, Mat}, $Q(F)$,
of $F$ is a quandle generated by $E=\{s_1,\dots,s_m\}$ 
with the following defining relations. 
Along a double point curve, let $s_j$ be the over-sheet
and $s_i$ (resp. $s_k$) the under-sheet which is behind 
(resp. in front of) the over-sheet $s_j$ with
respect to the normal vector of $s_j$. 
The defining relation is given by $s_i*s_j=s_k$ along 
the double point curve.
%(See Figure~\ref{}).
We note that $Q(F)$ is independent of the choice of the diagram
of $F$.
The following lemma
will be used to construct surface-knots satisfying 
the condition (i).

\begin{lem}\label{lem-quandle}
For surface-knots $F_0$ and $F$,
consider the connected sums $F_0 \# F$ and $F_0 \# -F^{*}$.
Then $Q(F_0 \# F)$ has the same presentation as $Q(F_0 \# -F^{*})$.
In particular, $Q(F_0 \# F)$ is isomorphic to $Q(F_0 \# -F^{*})$.
\end{lem}

\begin{proof}
A presentation of $Q(F_0 \# F)$ can be obtained from
that of $Q(F_0)$ and that of $Q(F)$
by adding a relation such as $a_0 = a$, where $a_0$ (resp. $a$)
is a generator of $Q(F_0)$ (resp. $Q(F)$)
corresponding to a sheet of a diagram of $F_0$ (resp. $F$). 
Since $Q(F)$ has the same presentation as $Q(-F^{*})$,
the result follows.
\end{proof}

\begin{rem}\label{rem-trefoil}\rm
The above lemma does not hold for classical knots in general.
Take the right-handed trefoils as $k_0$ and $k$ for example.
Then it is known in \cite[p.220]{Rolfsen} 
that the granny knot is not equivalent to the square knot up to orientation. 
(See Remark~\ref{rem-trefoil2} for an alternative proof of this fact.)
Thus we have that $Q(k_0 \# k)$ is not isomorphic to $Q(k_0 \# -k^{\ast})$.
\end{rem}

%%%%%%%%%%%%%%%%%%
\subsection{Quandle homology theory}
Before defining the fundamental class,
we briefly review the quandle homology theory 
defined in \cite{CJKLS}.
For $n>0$, let $\Hom{C}{R}_n(X)$ be the free 
abelian group generated by $n$-tuples
$(x_1,x_2,\dots,x_n)$ of elements of a quandle $X$. 
Put $\Hom{C}{R}_n(X)=0$ for $n\leq 0$. 
We define the boundary map $\Hom{\bou}{}_n : \Hom{C}{R}_{n}(X) 
\rightarrow \Hom{C}{R}_{n-1}(X)$ by
$$\begin{array}{llll}
\Hom{\bou}{}_n (x_1,\dots,x_{n})
& = & (-1)^{n-1} \sum_{i=1}^{n} (-1)^i & 
\bigl\{
(x_1,\dots, x_{i-1}, x_{i+1}, \dots,x_{n}) \\
&   &   &
-(x_1*x_i,\dots,x_{i-1}*x_i,x_{i+1},\dots,x_{n}) \bigr\} 
\end{array}$$
for $n>1$, and $\Hom{\bou}{}_n=0$ for $n\leq 1$.
It is easily verified that $\Hom{C}{R}_{*}(X)=
(\Hom{C}{R}_{n}(X),\Hom{\bou}{}_{n})$ is a chain complex.

For $n>1$, let $\Hom{C}{D}_n(X)$ be the submodule of 
$\Hom{C}{R}_n(X)$ generated by $n$-tuples $(x_1,x_2,\dots,x_n)$ 
with $x_i=x_{i+1}$ for some $i$ $(i=1,2,\dots, n-1)$.
Put $\Hom{C}{D}_n(X)=0$ for $n\leq 1$.
Quandle axiom (Q1) ensures that $\Hom{\bou}{}_n(\Hom{C}{D}_n(X)) 
\subset \Hom{C}{D}_{n-1}(X)$,
hence $\Hom{C}{D}_{*}(X)=(\Hom{C}{D}_{n}(X),\Hom{\bou}{}_{n})$ 
is a subcomplex of $\Hom{C}{R}_{*}(X)$.

Put $\Hom{C}{Q}_n(X) = \Hom{C}{R}_n(X) / \Hom{C}{D}_n(X)$ and 
$\Hom{C}{Q}_{*}(X)=(\Hom{C}{Q}_n(X), \Hom{\bou}{}_n)$,
where all the induced boundary operators are again denoted by 
$\Hom{\bou}{}_n$. 
For an element $x$ of $\Hom{C}{R}_n(X)$, 
we denote the equivalence class of $x$ 
by $x |_{\mathrm{Q}} \in \Hom{C}{Q}_n(X)$. 
The $n$th groups of cycles and boundaries of $\Hom{C}{Q}_{*}(X)$
are denoted by $\Hom{Z}{Q}_{n}(X)$ and $\Hom{B}{Q}_{n}(X)$, and
the $n$th homology group of this complex is
called the $n$th {\it quandle homology group} \cite{CJKLS}
and is denoted by $\Hom{H}{Q}_n(X)$.
For an abelian group $A$, define the cochain complex
\[ \begin{array}{ll}
\Coh{C}{W}^* (X;A) = \mbox{Hom}_{\Z}(\Hom{C}{W}_{*}(X),A) ,&
\Coh{\cob}{}^* = \mbox{Hom}(\Hom{\bou}{}_* , \mathrm{id}) 
\end{array} \]
in the usual way, where W $=$ R, D or Q.
The $n$th groups of cocycles and coboundaries of $\Coh{C}{Q}^{*}(X;A)$
are denoted by $\Coh{Z}{Q}^n(X;A)$ and $\Coh{B}{Q}^n(X;A)$, and
the $n$th cohomology group of this complex is
called the $n$th {\it quandle cohomology group} \cite{CJKLS}
and is denoted by $\Coh{H}{Q}^n(X;A)$.

%%%%%%%%%%%%%%%%%%%%%%%%%%%%%
\subsection{Fundamental classes}
Let $D$ be a diagram of a surface-link $F$
and let $E=\{s_1,\dots,s_m\}$ be the set of the sheets of $D$.
We often regard an element of $E$ as the element of the knot quandle $Q(F)$.

At a triple point $r$ of $D$, let $\vec{v_t}$,
$\vec{v_m}$ and $\vec{v_b}$ be the normal vectors to the top,
middle, and bottom sheet respectively. 
For the triple point $r$,
the {\it sign} $\varepsilon (r)$ is $1$ if the ordered triple 
$(\vec{v_t},\vec{v_m},\vec{v_b})$ matches the orientation of $\R^3$,
and $-1$ otherwise.

For a triple point $r$ of $D$, $C(r) = (s_b,s_m,s_t)$ is a triplet 
of elements of $Q(F)$, where $s_b$ is one of the four bottom sheets
from which the normal vectors of the middle and top sheets
point, $s_m$ is one of the two middle sheets from which 
the normal vector of the top sheet points, and $s_t$ is
the top sheet.

For a triple point $r$,
the {\it Boltzmann weight} $B(r) \in \Hom{C}{R}_3(Q(F))$ 
is defined by
$$B(r) :=  \varepsilon(r) C(r) \
\bigg( = \pm  (s_b,s_m,s_t) \bigg).$$
Let $|D| \in \Hom{C}{R}_3(Q(F))$ be the sum of the Boltzmann weights 
$B(r)$ of all triple points of the diagram $D$.
Then we have the following
(cf. \cite[Theorem $5.6$]{CJKLS}):
\begin{itemize}
\item
$|D|\big|_{\mathrm{Q}} \in \Hom{Z}{Q}_3(Q(F))$, and
\item
$|D'|\big|_{\mathrm{Q}}
-|D|\big|_{\mathrm{Q}} \in \Hom{B}{Q}_3(Q(F))$, 
for any other diagram $D'$ of $F$.
\end{itemize}
Thus the homology class of $|D|\big|_{\mathrm{Q}}$ 
is independent of the choice of the diagram $D$,
and the {\it fundamental class} \cite{CKS} (cf. \cite{Tana2}),
$[F]$, of a surface-link $F$ is defined by 
$$[F] := \Bigl[ |D|\big|_{\mathrm{Q}} \Bigr] \in \Hom{H}{Q}_3(Q(F)).$$

%%%%%%%%%%%%%%%%%%%%%%%%%%%
\subsection{Quandle cocycle invariants}\label{subsec-qc_inv}
Although a quandle cocycle invariant \cite{CJKLS} 
was originally introduced as an invariant for a surface-knot,
we use it as a tool for distinguishing given fundamental classes
(See Lemma~\ref{lem-fund} and Corollary~\ref{cor-fund} below).

Let $F$ be a surface-link and let $[F] \in \Hom{H}{Q}_3(Q(F))$ 
be the fundamental class of $F$.
For a finite quandle $X$, a abelian group $A$ and 
a $3$-cocycle $\theta \in \Coh{Z}{Q}^3(X;A)$,
we define a {\it quandle cocycle invariant} \cite{CJKLS}, 
$\Phi_\theta (F)$, by
$$\Phi_\theta (F) = \sum_{c:Q(F) \rightarrow X}
\langle\, {c}_*([F]), [\theta] \,\rangle \  \in \Z[A],$$
where ${c}_* : \Hom{H}{Q}_3(Q(F)) \rightarrow \Hom{H}{Q}_3(X)$
is a map induced from a quandle 
homomorphism $c:Q(F) \rightarrow X$, 
the element $[\theta]$ is a cohomology class of $\theta$, and
$$\langle \ , \ \rangle : \Hom{H}{Q}_3(X) \mathop{\otimes}\limits_{\Z} 
\Coh{H}{Q}^3(X;A) \rightarrow A$$ is a Kronecker product.
We note that the above summation is finite,
since the cardinarity of $X$ is finite.
The following are easy consequences of the construction of 
quandle cocycle invariants, and
Corollary~\ref{cor-fund} plays an important role 
in the proof of Theorem~\ref{th-quandle}.

\begin{lem}\label{lem-fund}
For surface-knots $F$ and $F'$,
if there exists a quandle isomorphism 
$f: Q(F) \rightarrow Q(F')$ such that $f_{*}[F]=[F']$,
then we have $\Phi_\theta(F)=\Phi_\theta(F')$ for
any finite quandle $X$, any abelian group $A$ 
and any $3$-cocycle $\theta$ of $Z_3^{\mathrm{Q}}(X;A)$.
\end{lem}

\begin{cor}\label{cor-fund}
For surface-knots $F$ and $F'$,
if there exists a finite quandle $X$,
an abelian group $A$ and a $3$-cocycle $\theta$
of $\Coh{Z}{Q}^3(X;A)$ such that
$$
\Phi_\theta(F)\neq \Phi_\theta(F')
\mbox{ and }
\Phi_\theta(F)\neq \Phi_\theta(-(F')^\ast ),
$$ 
then $F$ and $F'$ do not satisfy the condition \rm{(ii')}.
\end{cor}

%%%%%%%%%%%%%%%%%%%%%%%
\section{Proof of Theorem~\ref{th-quandle}}\label{sec-quandle}

Before proving Theorem~\ref{th-quandle},
we define two $S^2$-knots $F_{p,1}$ and
$F_{p,2}$, and study their properties. 
For an odd prime integer $p$,
let $K_p$ be the $2$-twist spun $S^2$-knot 
obtained from a $(2,p)$-torus knot.
Let $F_{p,1}$ be the connected sum of two copies of $K_p$,
and $F_{p,2}$ be the connected sum of $K_p$ and $-(K_p)^{\ast}$.

For a surface-knot $F$, let $\Phi_p(F)$ denote 
the quandle cocycle invariant of $F$
associated with Mochizuki's $3$-cocycle \cite{Mochi},
$\theta_p \in \Coh{Z}{Q}^3(R_p;\Z_p)$,
of the dihedral quandle $R_p$ 
and the coefficient group $\Z_p$. 
We note that the invariant $\Phi_p(F)$
takes values in $\Z[t,t^{-1}]/(t^p-1)$
$(\cong \Z[\Z_p])$.
Using Asami and Satoh's computation \cite{AS},
we have the following:
$$\Phi_p(F_{p,1}) = p \left(\sum_{k=0}^{p-1} t^{2k^2} \right)^2
\ \mbox{and} \quad
\Phi_p(F_{p,2}) = p \left(\sum_{k=0}^{p-1} t^{2k^2} \right) 
\left(\sum_{k=0}^{p-1} t^{-2k^2} \right).$$ 

\begin{prop}\label{prop-AS}
If p is an odd prime integer with $p \equiv 3 \pmod 4$, then
$\Phi_p(F_{p,1})$ is not equal to 
$\Phi_p(F_{p,2})$ in $\Z[t,t^{-1}]/(t^p-1)$.
\end{prop}

\begin{proof}
To compare thier values in $\Z[t,t^{-1}]/(t^p-1)$,
it is sufficient to  calculate \lq\lq constant terms\rq\rq\ 
of them, 
%where the \lq\lq constant term\rq\rq\ of 
%$\sum_{i} a_i t^i$ is defined by
%$$\sum_{i\equiv 0 \pmod p} a_i.$$
where the constant term of 
$\sum_{i} a_i t^i$  %$\in \Z[t,t^{-1}]/(t^p-1)$
is defined to be
$$\sum_{i\equiv 0 \pmod p} a_i  \ \  \in \Z .$$

For integers $i,j \in \{0,\dots,p-1\}$,  
it follows from the condition $p \equiv 3 \pmod 4$
that $2(i^2 + j^2) \equiv 0 \pmod p$ if and only if 
$(i,j)=(0,0)$.
Hence the constant term of $\Phi_p(F_{p,1})$ 
in $\Z[t,t^{-1}]/(t^p-1)$ is equal to $p$.

For integers $i,j \in \{0,\dots,p-1\}$,  
it is easy to see that
$2(i^2 - j^2) \equiv 0 \pmod p$ if and only if 
$$\begin{array}{lll}
(i,j)&=&(0,0),(1,1),\dots,(p-1,p-1), \\
&&(1,p-1),(2, p-2),\dots,(p-1,1).
\end{array}$$
Hence the constant term of $\Phi_p(F_{p,2})$ in 
$\Z[t,t^{-1}]/(t^p-1)$ is equal to $p(2p-1)$.
%Thus we conclude that $\Phi_p(F_{p,1})$ 
%is not equal to $\Phi_p(F_{p,2})$ in $\Z[t,t^{-1}]/(t^p-1)$.
\end{proof}

\begin{proof}[Proof of Theorem~\ref{th-quandle}]
We construct $S^2$-knots satisfying the condition
of Theorem~\ref{th-quandle}.
Let $\mathcal{P}$ be the set 
of odd prime integers with $p \equiv 3 \pmod 4$,
and take a subset $\{p_1,\dots ,p_n\}$ of $\mathcal{P}$
for any non-negative integer $n$.
We notice that the cardinality of $\mathcal{P}$ is countable.
Given an $n$-tuple $I = (e_1,\dots ,e_n) \in \{1,2\}^n$,
we consider the $S^2$-knot
$$F_I = F_{p_1,e_1} \# \dots \# F_{p_n,e_n},$$
and claim that these $2^n$ surface-knots
%if $I\neq I'$ then $F_I$ and $F_{I'}$ 
satisfy the condition.
For any two distinct elements $I=(e_1,\dots ,e_n)$ and 
$I'= (e'_1,\dots ,e'_n)$ of $\{1,2\}^n$,
we have $Q(F_I) \cong Q(F_{I'})$ 
by Lemma~\ref{lem-quandle},
that is, $F_I$ and $F_{I'}$ satisfy the condition (i).
Since $I\neq I'$,
there exists some $j$ $(j=1,\dots ,n)$ such that $e_j \neq {e'}_j$.
Thus we have 
$$\Phi_{p_j}(F_I) = \Phi_{p_j}(F_{p_j,e_j}) \neq 
\Phi_{p_j}(F_{p_j,{e'}_j}) = \Phi_{p_j}(F_{I'})$$
by Proposition~\ref{prop-AS}.
We can also show 
$$\Phi_{p_j}(F_I) \neq \Phi_{p_j}(-(F_{I'})^\ast )$$
in a similar way. 
Hence $F_I$ and $F_{I'}$ do
not satisfy the condition (ii')
by Corollary~\ref{cor-fund}.

When the genus $g$ is greater than zero, 
we consider the connected sum of $F_I$ 
and a trivial surface-knot of genus $g$.
Then these $2^n$ surface-knots of genus $g$
satisfy the condition of Theorem~\ref{th-quandle}.
\end{proof}

\begin{rem}\label{rem-trefoil2}\rm
We give an alternative proof of the fact mentioned 
in Remark~\ref{rem-trefoil}.
By the above proof, $F_{3,1}$ $(=K_3 \# K_3)$ 
is not equivalent to $F_{3,2}$ $(=K_3 \# -(K_3)^{\ast})$.
Then, for the right-handed trefoil knot 
(i.e., $(2,3)$-torus knot) $k_3$,
it follows from \cite{Lith} that $k_3 \# k_3$ is not equivalent to 
$k_3 \# -(k_3)^{\ast}$.
Since the trefoil knot is invertible,
the granny knot, $k_3 \# k_3$, is not equivalent to 
the square knot, $k_3 \# (k_3)^{\ast}$, up to orientation.
\end{rem}

%%%%%%%%%%%%%%%%%%%%%%%
\section{Proof of Theorem~\ref{th-fund}}\label{sec-fund}

The proof is divided into two cases:
One is the case where $g=0$ and 
the other is the case where $g>0$.

\subsection{$g=0$ case} 
Take integers $n,p,q>5$ such that
$p$ and $q$ are relatively prime.
Let $K$ be a $n$-twist spun $S^2$-knot 
obtained from a $(p,q)$-torus knot,
and $\widehat{K}$ be an $S^2$-knot obtained from
$K$ by Gluck surgery \cite{Glu}.
We remark that the exterior $E(K)$
of the $S^2$-knot $K$ is
homeomorphic to the exterior $E(\widehat{K})$ of 
$\widehat{K}$.
It is known in \cite{Gor} that 
the ambient space of $\widehat{K}$ is homeomorphic to 
the $4$-sphere $S^4$ and 
that $\widehat{K}$ is not equivalent to $K$
up to orientation.
In particular, $K$ and $\widehat{K}$
does not satisfy the condition (iii').

Let $\Sigma$ be the trivial surface-knot of genus two,
and consider the two surface-knots $K \# \Sigma$ and 
$\widehat{K} \# \Sigma$.
We notice that the exterior $E(K \# \Sigma)$ 
is homeomorphic to $E(\widehat{K} \# \Sigma)$.
Then $\widehat{K} \# \Sigma$ is equivalent to $K \# \Sigma$,
since a surface-knot of genus greater than one is determined by its
exterior \cite{Hi-Ka}.
Hence we have
$$
\begin{CD}
Q(K) @>\phi_1>\cong> Q(K\# \Sigma)  @>\phi_2>\cong> 
Q(\widehat{K} \# \Sigma) @>\phi_3>\cong> Q(\widehat{K})
\end{CD}
$$
and
$$
{(\phi_3 \circ \phi_2 \circ \phi_1)}_{\ast}[K]
=
{(\phi_3 \circ \phi_2)}_{\ast}[K\# \Sigma]
=
{(\phi_3)}_{\ast}[\widehat{K} \# \Sigma]
=[\widehat{K}],
$$
where the map $\phi_1$ (resp. $\phi_3$)
is induced by doing the connected sum of the trivial surface-knot 
$\Sigma$ to $K$ (resp. $\widehat{K}$),
and the map $\phi_2$ is induced from
the equivalence between $K \# \Sigma$
and $\widehat{K} \# \Sigma$.
When we vary integers $n$, $p$ and $q$,
we can obtain infinitely many such pairs.

%%%%%%%%%%%%%%%%%%%%%%%%%%
\subsection{$g>0$ case}\label{subsec-}
Let $T(k)$ denote the spun $T^2$-knot obtained from 
a non-trivial classical knot $k$, and 
let $\tilde{T}(k)$ denote the turned spun $T^2$-knot 
%\cite{}
obtained from $k$.
Take a ribbon surface-knot $G$ of genus $g-1$ ($\geq 0$)
and consider the two surface-knots 
$G \# T(k)$ and $G \# \tilde{T}(k)$ of genus $g$.
It is easy to see that these two surface-knots
satisfy the condition (ii).
We note that the fundamental classes of them are
equal to zero elements.

To distinguish them,
we use Kawauchi's Gauss sum invariant \cite[p.1047]{Kawa},
$\varsigma(F) \in \Z$, of a surface-knot $F$.
It is known in \cite{Kawa} that $\varsigma(G)=2^{g-1}$,
$\varsigma(T(k))=2$ and $\varsigma(\tilde{T}(k))=0$. 
Using the connected sum formula \cite[Theorem $1.2$]{Kawa} 
$$\varsigma(F_1 \# F_2)=\varsigma(F_1)\varsigma(F_2),$$
we have
$$\varsigma(G \# T(k))=2^g \neq 
0 =\varsigma(G \# \tilde{T}(k)),$$ 
and it follows that they do not satisfy the condition (iii').
When we vary a non-trivial classical knot $k$,
we can obtain infinitely many such pairs.

\begin{rem}\rm
We may take any surface-knot $G$ of genus $g-1$
as long as it satisfies the condition $\varsigma(G)\neq 0$,
though we take a ribbon surface-knot as $G$ 
in the above proof.
\end{rem}

%%%%%%%%%%%%%%%%%%%
% Acknowledgments
%%%%%%%%%%%%%%%%%%
\section*{Acknowledgments}
The author would like to express his sincere gratitude
to Yukio Matsumoto for encouraging him.
He would like to thank Seiichi Kamada
for helpful comments on Theorem~\ref{th-quandle}, 
Akira Yasuhara for telling me the result \cite{Hi-Ka}
due to Jonathan A. Hillman and Akio Kawauchi,
and Isao Hasegawa for stimulating discussions.
He would also like to thank Eri Hatakenaka for her natural
question to him: How powerful is the knot quandle 
for surface-knots?
This research is supported by JSPS Research Fellowships for
Young Scientists.

%%%%%%%%%%%%%%%%
% bibliography
%%%%%%%%%%%%%%%

%\bibliographystyle{amsplain}

%%%%%%%%%%%%%%%%%%%%%%%%%%%%%%%%

\begin{thebibliography}{99}

\bibitem{AS}
S. Asami and S. Satoh,
\textit{An infinite family of non-invertible surfaces in 4-space,}
Bull. London Math. Soc. \textbf{37} (2005), no. 2, 285--296. 

\bibitem{BL}
J. Bojarczuk and P. Lopes, 
\textit{Quandles at finite temperatures. III,} 
J. Knot Theory Ramifications \textbf{14} (2005), no. 3, 275--373.

\bibitem{CJKLS}
J. S. Carter, D. Jelsovsky, S. Kamada, 
L. Langford and M. Saito, 
\textit{Quandle cohomology and state-sum invariants 
of knotted curves and surfaces,}
Trans. Amer. Math. Soc. \textbf{355} (2003), no. 10, 3947--3989.

\bibitem{CJKS}
J. S. Carter, D. Jelsovsky, S. Kamada and M. Saito, 
\textit{Computations of quandle cocycle invariants of knotted curves and surfaces,} 
Adv. Math. \textbf{157} (2001), no. 1, 36--94. 

\bibitem{CKS}
J. S. Carter, S. Kamada, M. Saito,
\textit{Diagrammatic computations for quandles and 
cocycle knot invariants,}
Contemp. Math., {\bf 318}, 51--74

\bibitem{CS-book} 
J. S. Carter and M. Saito, 
``Knotted surfaces and their diagrams'',
Math. Surveys and Monographs \textbf{55}, 
Amer. Math. Soc., 1998.

\bibitem{CSS}
J. S. Carter, M. Saito, and S. Satoh, 
\textit{Ribbon concordance of surface-knots via quandle cocycle invariants,}
to appear in J. Aust. Math. Soc.


\bibitem{Eis}
M. Eisermann, 
\textit{Homological characterization of the unknot,}
J. Pure Appl. Algebra \textbf{177} (2003), no. 2, 131--157. 

\bibitem{Glu}
H. Gluck,
\textit{The embedding of two-spheres in the four-sphere,} 
Trans. Amer. Math. Soc. \textbf{104} (1962) 308--333.

\bibitem{Gor}
C. McA. Gordon, 
\textit{Knots in the $4$-sphere,}
Comment. Math. Helv. \textbf{51} (1976), no. 4, 585--596.

\bibitem{Hat}
E. Hatakenaka,
\textit{An estimate of the triple point numbers of surface-knots 
by quandle cocycle invariants,}
Topology Appl. \textbf{139} (2004), no. 1-3, 129--144.

\bibitem{Hi-Ka}
J. A. Hillman and A. Kawauchi,
\textit{Unknotting orientable surfaces in the $4$-sphere,} 
J. Knot Theory Ramifications \textbf{4} (1995), no. 2, 213--224.

\bibitem{Iwa}
M. Iwakiri,
\textit{Calculation of dihedral quandle cocycle invariants 
of twist spun 2-bridge knots,}
J. Knot Theory Ramifications \textbf{14} (2005), no. 2, 217--229.

\bibitem{Iwa2}
M. Iwakiri,
\textit{Unknotting and triple point cancelling numbers of surface links,}
preprint.

\bibitem{Joy}
D. Joyce, 
\textit{A classifying invariant of knots, the knot quandle,}
J. Pure Appl. Algebra \textbf{23} (1982), no. 1, 37--65.

\bibitem{Kawa}
A. Kawauchi,
\textit{On pseudo-ribbon surface-links,} 
J. Knot Theory Ramifications \textbf{11} (2002), no. 7, 1043--1062. 

\bibitem{Lith}
R. A. Litherland, 
\textit{Symmetries of twist-spun knots,}
Knot theory and manifolds (Vancouver, B.C., 1983), 97--107, 
Lecture Notes in Math., \textbf{1144}, Springer, Berlin, 1985.


\bibitem{Mat}
S. V. Matveev, 
\textit{Distributive groupoids in knot theory,}
(Russian) Mat. Sb. (N.S.) \textbf{119(161)} 
(1982), no. 1, 78--88, 160.

\bibitem{Mochi}
T. Mochizuki,
\textit{Some calculations of cohomology groups of finite Alexander quandles,}
J. Pure Appl. Algebra \textbf{179} (2003), no. 3, 287--330. 

\bibitem{Rolfsen}
D. Rolfsen,
\textit{Knots and links,} 
Mathematics Lecture Series, No. 7. Publish or Perish, Inc., 
Berkeley, Calif., 1976


\bibitem{Sai-Sat}
M. Saito and S. Satoh,
\textit{The Spun Trefoil Needs Four Broken Sheets,}
to appear in J. Knot Theory Ramifications, 
\textbf{14} (2005), no. 7. 

\bibitem{Sa-Shi}
S. Satoh and A. Shima,
\textit{The 2-twist-spun trefoil has the triple point number four,}
Trans. Amer. Math. Soc. \textbf{356} (2004), no. 3, 1007--1024.

\bibitem{Sa-Shi2}
S. Satoh and A. Shima,
\textit{Triple point numbers and quandle cocycle invariants 
of knotted surfaces in 4-space,}
New Zealand J. Math. \textbf{34} (2005), no. 1, 71--79.

\bibitem{Suc}
A. I. Suciu, 
\textit{Infinitely many ribbon knots with the same fundamental group,} 
Math. Proc. Cambridge Philos. Soc. \textbf{98} (1985), no. 3, 481--492.

\bibitem{Tana}
K. Tanaka,
\textit{The braid index of surface-knots and quandle colorings,} 
Illinois J. Math, \textbf{49} (2005), no. 2, 517--522. 

\bibitem{Tana2}
K. Tanaka,
\textit{On surface-links represented by diagrams 
with two or three triple points,} 
to appear in J. Knot Theory Ramifications, 
\textbf{14} (2005), no. 8. 


\end{thebibliography}
\end{document}